\documentclass[12pt,reqno]{amsart}
\usepackage{enumerate, latexsym, amsmath, amsfonts, amssymb, amsthm, color}
\def\pmod #1{\ ({\rm{mod}}\ #1)}

\def\Q{\Bbb Q}

\def\l{\left}
\def\r{\right}
\def\bg{\bigg}
\def\({\bg(}
\def\){\bg)}
\def\t{\text}
\def\f{\frac}

\def\ls{\leqslant}

\def\ve{\varepsilon}

\def\eq{\equiv}

\def\Proof{\noindent{\it Proof}}

\theoremstyle{plain}
\newtheorem{theorem}{Theorem}

\newtheorem{lemma}{Lemma}

\newtheorem{conjecture}{Conjecture}
\theoremstyle{definition}

\theoremstyle{remark}
\newtheorem{remark}{Remark}

\makeatletter
\@namedef{subjclassname@2020}{%
  \textup{2020} Mathematics Subject Classification}
\makeatother
 \vspace{4mm}

\begin{document}

\medskip

\title
[Evaluation of a determinant involving Legendre symbols]
{Evaluation of a determinant involving Legendre symbols}

\author
[Chen-Kai Ren and Zhi-Wei Sun] {Chen-Kai Ren and Zhi-Wei Sun}

\address {(Chen-Kai Ren) School of Mathematics, Nanjing
University, Nanjing 210093, People's Republic of China}
\email{ckren@smail.nju.edu.cn}

\address{(Zhi-Wei Sun, corresponding author) School of Mathematics, Nanjing
University, Nanjing 210093, People's Republic of China}
\email{zwsun@nju.edu.cn}

\keywords{Determinant, Legendre symbol, quadratic Gauss sum, Lagrange's interpolation formula.
\newline \indent 2020 {\it Mathematics Subject Classification}. Primary 11C20, 15A15; Secondary 11A15.
\newline \indent Supported by the Natural Science Foundation of China (grant 12371004).}

\begin{abstract} Let $p>3$ be a prime, and let $(\frac{\cdot}p)$ be the Legendre symbol.
Let $A_p(x)$ denote the matrix $[x+a_{ij}]_{1\leqslant i,j\leqslant (p-1)/2}$, where
$$  a_{ij}=\begin{cases}
    (\frac{j}{p}) &\text{if} \ i=1, \\
    (\frac{i+j}{p}) &\text{if} \ i>1.
\end{cases}
$$
In 2018 Z.-W. Sun conjectured that $\det A_p(0)=-2^{(p-3)/2}$ if $p\equiv 3 \pmod{4}$, which was later confirmed by G. Zaimi. In this paper we evaluate $\det A_p(x)$ completely.
\end{abstract}
\maketitle

\section{Introduction}
\setcounter{lemma}{0}
\setcounter{theorem}{0}
\setcounter{corollary}{0}
\setcounter{remark}{0}
\setcounter{equation}{0}

For an $n\times n$ matrix $M=[a_{ij}]_{1\leqslant i,j\leqslant n}$ over a field, as usual $\det M$  denotes the determinant of $M.$

Let $p$ be an odd prime. Let $\ve_p$ and $h(p)$ denote the fundamental unit and the class number
of the real quadratic field $\Q(\sqrt p)$, respectively. Write
\begin{equation}\label{vp}
   \ve_p^{h(p)}=a_p+b_p\sqrt p,
\end{equation}
where $a_p,b_p\in \mathbb{Q}.$
Let $(\frac{.}{p})$ be the Legendre symbol.
In 2004, via quadratic Gauss sums, R. Chapman \cite{C1} showed that if $p>3$ then
$$\det\l[\l(\f{i+j-1}p\r)\r]_{1\ls i,j\ls(p-1)/2}=\begin{cases}
(-1)^{(p-1)/4}2^{(p-1)/2}b_p&\t{if}\ p\eq1\pmod4,
\\0&\t{if}\ p\eq3\pmod4,\end{cases}$$
and
$$\det\l[\l(\f{i+j-1}p\r)\r]_{1\ls i,j\ls(p+1)/2}=\begin{cases}
(-1)^{(p+3)/4}2^{(p-1)/2}a_p&\t{if}\ p\eq1\pmod4,
\\2^{(p-1)/2}&\t{if}\ p\eq3\pmod4.\end{cases}$$
As $(\f{p+1}2-i)+(\f{p+1}2-j)-1\eq-(i+j)\pmod p$, we have
$$\det\bg[\l(\f{i+j-1}p\r)\bg]_{1\ls i,j\ls(p-1)/2}=\l(\f{-1}p\r)\det\bg[\l(\f{i+j}p\r)\bg]_{1\ls i,j\ls(p-1)/2}$$
and
$$\det\bg[\l(\f{i+j-1}p\r)\bg]_{1\ls i,j\ls(p+1)/2}
=\det\bg[\l(\f{i+j}p\r)\bg]_{0\ls i,j\ls(p-1)/2}.$$
Thus Chapman \cite{C1} actually determined the values of
$$\det\bg[\l(\f{i+j}p\r)\bg]_{1\ls i,j\ls(p-1)/2}
\ \ \t{and}\ \ \det\bg[\l(\f{i+j}p\r)\bg]_{0\ls i,j\ls(p-1)/2}.$$

 Let $p>3$ be a prime. Motivated by Chapamn's work \cite{C1},  in 2018 Z.-W. Sun \cite{SMa} considered
 the matrix $A_p$ which is obtained by replacing the first row in $[(\f{i+j}p)]_{1\ls i,j\ls(p-1)/2}$
by $((\f1p),\ldots,(\f{(p-1)/2}p))$. He conjectured that $\det A_p=-2^{(p-3)/2}$ if $p\equiv 3 \pmod{4}$, and this was later confirmed by  G. Zaimi (cf. the answer in \cite{SMa}).

In this paper we obtain the following general result.

\begin{theorem} \label{Th1.1}
    Let $p>3$ be a prime, and  let $A_p(x)=[x+a_{ij}]_{1\leqslant i,j\leqslant (p-1)/2},$
where
$$a_{ij}=\begin{cases}
    (\frac{j}{p}) &\t{if} \ i=1, \\
    (\frac{i+j}{p}) &\t{if} \ i>1.
\end{cases}
$$
Then
$$\det A_p(x)=\begin{cases}
    2^{(p-3)/2}(\frac{2}{p})(a_p-b_p+(a_p-pb_p)x) & \t{if} \ p\equiv 1 \pmod{4}, \\
  2^{(p-3)/2}(x-1)   & \t{if} \ p\equiv 3 \pmod{4},
\end{cases}
$$
where $a_p$ and $b_p$ are given by \eqref{vp}.
\end{theorem}
\begin{remark} Let $p>3$ be a prime. For the matrix $M_p^+$ obtained from
$[(\f{i+j}p)]_{0\ls i,j\ls(p-1)/2}$ via replacing all the entries in the first row by $1$.
Sun \cite[Conjecture 4.6]{S19} conjectured that
$$\det M_p^+=\begin{cases}(-1)^{(p-1)/4}2^{(p-1)/2}&\t{if}\ p\eq1\pmod4,
\\(-1)^{(h(-p)-1)/2}2^{(p-1)/2}&\t{if}\ p\eq3\pmod4,
\end{cases}$$
where $h(-p)$ denotes the class number of the imaginary quadratic field $\Q(\sqrt{-p})$.
This was confirmed by L.-Y. Wang and H.-L. Wu \cite[Theorem 1.1]{W22}.
\end{remark}

However, our method to prove Theorem \ref{Th1.1} does not work for the following conjecture
made by the second author.

\begin{conjecture} [Sun, 2024]
    Let $p>3$ be a prime, and let
    $C_p(x)=[x+c_{ij}]_{1\leqslant i,j\leqslant (p-1)/2},$
    where
    $$  c_{ij}=\begin{cases}
    (\frac{j}{p}) &\t{if} \ i=1, \\
    (\frac{i-j}{p}) &\t{if} \ i>1.
\end{cases}
$$
 If $p\eq1\pmod4$, then
 $$\det C_p(x)=\l(\f 2p\r)\l(a_p'-\frac{p+1}{2}b_p'\r)+d_px$$
 for some integer $d_p$, 
 where we write $\ve_p^{(2-(\f 2p))h(p)}$ as $a_p'+b_p'\sqrt p$ with $a_p',b_p'\in\Q$.
 When $p\eq3\pmod4$, we have
    $$\det C_p(x)=(-1)^{(p-3)/4}+\l(2\l(\frac{2}{p}\r)(-1)^{(h(-p)-1)/2}-\frac{p-1}{2}\r)x.$$
\end{conjecture}
\begin{remark} This conjecture is motivated by the following result conjectured by Chapman and proved by  M. Vsemirnov \cite{V12,V13}: For any odd prime $p$, we have
$$\det\left[\left(\frac{i-j}{p}\right)\right]_{0\leqslant i,j\leqslant (p-1)/{2}}
=\begin{cases}-a_p'&\t{if}\ p\eq1\pmod 4,\\1&\t{if}\ p\eq3\pmod 4,
\end{cases}$$
For any prime $p>3$, let $M_p$ be the matrix obtained from  $[(\f{i-j}p)]_{0\ls i,j\ls (p-1)/2}$
via replacing all entries in the first row by $1$. Sun \cite{S19} conjectured that
$$\det M_p=\begin{cases}(-1)^{(p-1)/4}&\t{if}\ p\eq1\pmod4,
\\(-1)^{(h(-p)-1)/2}&\t{if}\ p\eq3\pmod4,\end{cases}$$
and this was confirmed by Wang and Wu \cite{W25}.
\end{remark}

 We are going to provide some auxiliary results in the next section and prove Theorem 1.1 in Section 3.

\section{Some auxiliary results}
\setcounter{lemma}{0}
\setcounter{theorem}{0}
\setcounter{corollary}{0}
\setcounter{remark}{0}
\setcounter{equation}{0}

The following result is well-known, and it can be found in \cite{V16}.

\begin{lemma}[Matrix-Determinant Lemma]
If $H$ is  an $n\times n$ matrix over the field $\mathbb{C}$ of complex numbers, and $\textbf{u}$ and $\textbf{v}$ are two $n$-dimensional column vectors with entries in $\mathbb{C}$, then
   $$\det(H+\textbf{u}\textbf{v}^T)=\det H +\textbf{v}^T\mathrm{adj}(H)\textbf{u}, $$
   where $\mathrm{adj}(H)$ denotes the adjugate matrix of $H$.
\end{lemma}

 We also need the following  auxiliary result on roots of unity from \cite[Theorem 1.3]{S20}.
\begin{lemma} Let $p$ be an odd prime, and set $\zeta=e^{2\pi i/p}$.

{\rm (i)} If $p\equiv 1 \pmod{4}$, then
\begin{equation}\label{z1}
    \prod_{k=1}^{(p-1)/2}(1-\zeta^{k^2})=\sqrt{p}\varepsilon_p^{-h(p)},
\end{equation}
and
   \begin{equation}\label{z2}
    \prod_{1\leqslant j<k\leqslant (p-1)/2 }(\zeta^{j^2}-\zeta^{k^2})^2=(-1)^{(p-1)/4}p^{(p-3)/4}\varepsilon_p^{h(p)}.
   \end{equation}

{\rm (ii)} When $p\equiv 3 \pmod{4}$, we have
 \begin{equation}\label{z3}
    \prod_{k=1}^{(p-1)/2}(1-\zeta^{k^2})=(-1)^{(h(-p)+1)/2}\sqrt{p}i,
     \end{equation}
 and    \begin{equation}\label{z4}
    \prod_{1\leqslant j<k\leqslant (p-1)/2 }(\zeta^{j^2}-\zeta^{k^2})^2=(-p)^{(p-3)/4},
   \end{equation}
   where $h(-p)$ is  the class number of the imaginary quadratic field $\mathbb{Q}(\sqrt{-p}).$
\end{lemma}

 The following result is classical.
 
\begin{lemma}[Lagrange's interpolation formula] Let $x_0,x_1,\cdots,x_{n},y_0,y_1,\cdots,y_{n}$ be complex numbers with all $x_0,x_1,\ldots,x_n$ distinct. Then there exists a unique polynomial $P_n(x)\in \mathbb{C}[x]$ with $\mathrm{deg} P_n(x)\leqslant n$ such that $P_n(x_i)=y_i$ for all $i=0,\ldots,n$. In fact, $P_n(x)$ is given by
$$P_n(x)=\sum_{k=0}^{n}y_k\cdot \prod_{\substack{j=0\\j\neq k}}^{n+1}\frac{x-x_j}{x_k-x_j}.$$
\end{lemma}

Let $p$ be an odd prime and set $\zeta=e^{2\pi i/p}$. For the quadratic Gauss sum
$$ \tau_p=\sum_{k=1}^{p-1}\l(\frac{k}{p}\r)\zeta^k.$$
It is well known that
$$ \tau_p=\begin{cases}
    \sqrt{p}  &\t{if}\ p\eq1\pmod 4,
\\ i\sqrt{p}  &\t{if}\ p\eq3\pmod 4.\\
\end{cases}$$
Moreover, for any integer $a$ not divisible by $p$, we have
\begin{equation}\label{ta}
  \l(\frac{a}{p}\r)\tau_p=\sum_{k=0}^{p-1}\zeta^{ak^2}=1+2\sum_{k=1}^{(p-1)/2}\zeta^{ak^2}.
\end{equation}

\begin{lemma}
    Let $p$ be a prime with $p\eq1\pmod4$, and set $\zeta=e^{2\pi i/p}$. Then
    \begin{equation}\label{l1}
     \sum_{k=1}^{(p-1)/2}\frac{1}{\zeta^{2k^2}} \prod_{\substack{j=1\\j\neq k}}^{(p-1)/2}\(\frac{1-\zeta^{j^2}}{\zeta^{k^2}-\zeta^{j^2}}\)^2= \varepsilon_p^{-2h(p)}-1 ,
     \end{equation}
   and  \begin{equation}\label{l2}
         \sum_{k=1}^{(p-1)/2}\frac{1}{\zeta^{3k^2}} \prod_{\substack{j=1\\j\neq k}}^{(p-1)/2}\(\frac{1-\zeta^{j^2}}{\zeta^{k^2}-\zeta^{j^2}}\)^2=(\sqrt{p}+1)\varepsilon_p^{-2h(p)}-1.
     \end{equation}
\end{lemma}
\Proof. For convenience, we set
$$  R=\l\{1\leqslant j\leqslant p-1: \(\frac{j}{p}\)=1\r\}$$
and
$$  N=\l\{1\leqslant j\leqslant p-1: \(\frac{j}{p}\)=-1\r\}.$$
Then
\begin{equation}\label{P}
  P(x):=\frac{1}{x}\(\prod_{j\in N}(x-\zeta^j)-\prod_{k\in R}(x-\zeta^k)\).
\end{equation}
 is a polynomial with $\mathrm{deg}  P(x)<(p-1)/2. $

We now prove \eqref{l1}. Applying Lemma 2.3 on $P(x)$ at the points $x=\zeta^{1^2},\cdots,\zeta^{(\frac{p-1}{2})^2},$  we obtain
\begin{equation}\label{gl1}
  P(x)=\sum_{k=1}^{(p-1)/2}P(\zeta^{k^2}) \prod_{\substack{j=1\\j\neq k}}^{(p-1)/2}\frac{x-\zeta^{j^2}}{\zeta^{k^2}-\zeta^{j^2}}.
\end{equation}
For each $k=1,\ldots,(p-1)/2$, clearly
$$P(\zeta^{k^2})=\frac{1}{\zeta^{k^2}}\prod_{j\in N}(\zeta^{k^2}-\zeta^j)$$
and
\begin{equation}
\begin{aligned}
 \prod_{\substack{j=1\\j\neq k}}^{(p-1)/2}(\zeta^{k^2}-\zeta^{j^2})&=\frac{1}{(\zeta^{k^2}-1)\prod_{j\in N}(\zeta^{k^2}-\zeta^j)} \prod_{\substack{j=0\\j\neq k^2}}^{p-1}(\zeta^{k^2}-\zeta^{j}) \\
 &=\frac{1}{(\zeta^{k^2}-1)\prod_{j\in N}(\zeta^{k^2}-\zeta^j)}\lim_{x\to\zeta^{k^2}}\frac{x^p-1}{x-\zeta^{k^2}}   \\
  &= \frac{-p}{\zeta^{k^2}(1-\zeta^{k^2})\prod_{j\in N}(\zeta^{k^2}-\zeta^j)}.
\end{aligned}
\end{equation}
Taking $x=1$ in \eqref{gl1}, we obtain
\begin{equation}
  \begin{split}
    P(1)&=\sum_{k=1}^{(p-1)/2}\frac{1}{\zeta^{k^2}}\prod_{j\in N}(\zeta^{k^2}-\zeta^j) \prod_{\substack{j=1\\j\neq k}}^{(p-1)/2}\frac{1-\zeta^{j^2}}{\zeta^{k^2}-\zeta^{j^2}}  \\
    &=\sum_{k=1}^{(p-1)/2}\frac{-p}{\zeta^{2k^2}(1-\zeta^{k^2})} \prod_{\substack{j=1\\j\neq k}}^{(p-1)/2}\frac{1-\zeta^{j^2}}{(\zeta^{k^2}-\zeta^{j^2})^2} \\
    &=\frac{-p}{\prod_{j=1}^{(p-1)/2}(1-\zeta^{j^2})}\sum_{k=1}^{(p-1)/2}\frac{1}{\zeta^{2k^2}} \prod_{\substack{j=1\\j\neq k}}^{(p-1)/2}\(\frac{1-\zeta^{j^2}}{\zeta^{k^2}-\zeta^{j^2}}\)^2.
  \end{split}
\end{equation}
Since
$$ \prod_{r=1}^{p-1}(1-\zeta^r)=\lim_{x\to1} \frac{x^p-1}{x-1}=p,
$$
from \eqref{z1} we get that
\begin{equation}\label{w1}
  \prod_{j\in N}(1-\zeta^j)= \varepsilon_p^{h(p)} \sqrt{p}.
\end{equation}
Combining Lemma 2.2, (2.8),
 (2.11) and \eqref{w1} we obtain
\begin{align*}
  &\sum_{k=1}^{(p-1)/2}\frac{1}{\zeta^{2k^2}} \prod_{\substack{j=1\\j\neq k}}^{(p-1)/2}\(\frac{1-\zeta^{j^2}}{\zeta^{k^2}-\zeta^{j^2}}\)^2 \\
 =&-\frac{1}{p}\(\prod_{j=1}^{(p-1)/2}(1-\zeta^{j^2})\prod_{k\in N}(1-\zeta^k)-\prod_{j=1}^{(p-1)/2}(1-\zeta^{j^2})^2\)  \\
 =& -\frac{1}{p}\(p-p\varepsilon_p^{-2h(p)}\)=\varepsilon_p^{-2h(p)}-1.
\end{align*}
Thus \eqref{l1} holds.

 We now turn to prove \eqref{l2}. Applying Lemma 2.3 on $P(x)$ at the points $x=0, \zeta^{1^2},\cdots,\zeta^{(\frac{p-1}{2})^2},$ we get
 $$
  P(x)=P(0)\prod_{k=1}^{(p-1)/2}\frac{x-\zeta^{k^2}}{0-\zeta^{k^2}}+\sum_{k=1}^{(p-1)/2} P(\zeta^{k^2})\frac{x}{\zeta^{k^2}}\prod_{\substack{j=1\\j\neq k}}^{(p-1)/2}\frac{x-\zeta^{j^2}}{\zeta^{k^2}-\zeta^{j^2}}.
$$
 Taking $x=1$ in the last equality, we obtain
\begin{equation}\label{p2}
    P(1)=P(0)\prod_{k=1}^{(p-1)/2}\frac{1-\zeta^{k^2}}{0-\zeta^{k^2}}+\sum_{k=1}^{(p-1)/2}\frac{1}{\zeta^{2k^2}}\prod_{j\in N}(\zeta^{k^2}-\zeta^j) \prod_{\substack{j=1\\j\neq k}}^{(p-1)/2}\frac{1-\zeta^{j^2}}{\zeta^{k^2}-\zeta^{j^2}}.
\end{equation}
In view of (2.10), we can verify that
\begin{equation}\label{p3}
   P(1)-P(0)\prod_{k=1}^{\frac{p-1}{2}}(1-\zeta^{k^2}) =\frac{-p}{\prod_{j=1}^{\frac{p-1}{2}}(1-\zeta^{j^2})} \sum_{k=1}^{\frac{p-1}{2}}\frac{1}{\zeta^{3k^2}} \prod_{\substack{j=1\\j\neq k}}^{\frac{p-1}{2}}\(\frac{1-\zeta^{j^2}}{\zeta^{k^2}-\zeta^{j^2}}\)^2.
\end{equation}

Now, we calculate $P(0)$. From the definition of $P(x)$, we obtain
\begin{align*}
    P(0)&= \(\prod_{j\in N}(x-\zeta^j)-\prod_{k\in R}(x-\zeta^k)\)^{'}\Big|_{x=0}\\
    &=-\sum_{j\in N}\zeta^{-j}\prod_{i\in N}\zeta^i+\sum_{k\in R}\zeta^{-k}\prod_{l\in R}\zeta^l.\\
\end{align*}
Since
$$\sum_{k=1}^{p-1}k=\frac{p(p-1)}{2} \equiv 0 \pmod{p} \ \t{and}\ \sum_{k=1}^{(p-1)/2}k^2=\frac{p(p^2-1)}{24}\equiv 0 \pmod{p},
$$
we have
\begin{equation}\label{pi}
   \prod_{i\in N}\zeta^i=\prod_{l\in R}\zeta^l=1,
\end{equation}
and hence
\begin{equation}
   \begin{aligned}
    P(0)&= -\sum_{j\in N}\zeta^{-j}+\sum_{k\in R}\zeta^{-k} \\
    &=-\sum_{j=1}^{p-1}\zeta^{-j}+2\sum_{k=1}^{(p-1)/2}\zeta^{-k^2}\\
    &=1+2\sum_{k=1}^{(p-1)/2}\zeta^{k^2}=\tau_p=\sqrt{p}.
\end{aligned}
\end{equation}
Combining Lemma 2.2, (2.8), (2.12), \eqref{p3} and (2.16), we obtain
\begin{align*}
  \sum_{k=1}^{(p-1)/2}\frac{1}{\zeta^{3k^2}} \prod_{\substack{j=1\\j\neq k}}^{(p-1)/2}\(\frac{1-\zeta^{j^2}}{\zeta^{k^2}-\zeta^{j^2}}\)^2  &=
\varepsilon_p^{-2h(p)}-1+\frac{\sqrt{p}}{p}\prod_{k=1}^{\frac{p-1}{2}}(1-\zeta^{k^2})^2\\
&=\varepsilon_p^{-2h(p)}-1+\frac{\sqrt{p}}{p}\l(\sqrt{p}\varepsilon_p^{-h(p)}\r)^2\\
&=(\sqrt{p}+1)\varepsilon_p^{-2h(p)}-1.
\end{align*}
This concludes  the proof of \eqref{l2}.
\qed

The following auxiliary result from \cite{SMa} deals with primes  $p\equiv 3\pmod{4}.$

\begin{lemma}[Zaimi]
    Let $p$ be a  prime with $p\equiv 3 \pmod{4}$ and let $\zeta=e^{2\pi i/p}$. Then
    \begin{equation}\label{l3}
     \sum_{k=1}^{(p-1)/2}\frac{1}{\zeta^{2k^2}}\cdot \prod_{\substack{j=1\\j\neq k}}^{(p-1)/2}\(\frac{1-\zeta^{j^2}}{\zeta^{k^2}-\zeta^{j^2}}\)^2= -2
     \end{equation}
     and
     \begin{equation}\label{l4}
         \sum_{k=1}^{(p-1)/2}\frac{1}{\zeta^{3k^2}}\cdot \prod_{\substack{j=1\\j\neq k}}^{(p-1)/2}\(\frac{1-\zeta^{j^2}}{\zeta^{k^2}-\zeta^{j^2}}\)^2=-2+i\sqrt{p}.
     \end{equation}
\end{lemma}

\section{proof of Theorem 1.1}
\setcounter{lemma}{0}
\setcounter{theorem}{0}
\setcounter{corollary}{0}
\setcounter{remark}{0}
\setcounter{equation}{0}
 
 Set $n=(p-1)/2.$ Clearly $p\nmid i+j$ for any $1\leqslant i,j\leqslant n.$ Let $\zeta=e^{2\pi i/p}$ and $ \tau_p=\sum_{k=1}^{p-1}(\frac{k}{p})\zeta^k$.
From \eqref{ta}, we have
$$ \(\frac{i+j}{p}\)=\frac{1}{\tau_p}\(1+2\sum_{k=1}^{n}\zeta^{ik^2}\cdot\zeta^{jk^2}\).
$$
Define $ B_p=[b_{ij}]_{1\leqslant i,j\leqslant n}$, where
$$b_{ij}=\begin{cases}
    x\tau_p+1+ 2\sum_{k=1}^{n}\zeta^{jk^2}&\t{if} \ i=1, \\
  x\tau_p+1+ 2\sum_{k=1}^{n}\zeta^{(i+j)k^2}   &\t{if} \ i>1. \\
\end{cases}
$$
Obviously,
\begin{equation}\label{bp'}
 \det B_p=\tau_p^{n}\det A_p(x).
\end{equation}
Let $\textbf{v}$ be the $n$-dimensional column vector with all entries 1. Then
$$B_p=(x\tau_p+1)\textbf{v}\textbf{v}^T+2MN^T,
$$
where
$$ N=[\zeta^{ij^2}]_{1\leqslant i,j\leqslant n}$$
and $M$ is the matrix obtaining from $N$ via replacing all the entries in the first row by 1. By Lemma 2.1, we have
\begin{equation}
    \begin{aligned}
    \det B_p&=\det(2MN^T)+(x\tau+1)\textbf{v}^T\mathrm{adj}(2MN^T)\textbf{v}\\
    &=2^n\det M\det N+2^{n-1}(x\tau+1)\textbf{v}^T\mathrm{adj}(N)^T\mathrm{adj}(M)\textbf{v}.
\end{aligned}
\end{equation}

First, we calculate $\det M$ and $\det N.$ For convenience, let us write $\eta_i=\zeta^{i^2}$ for any $1\leqslant i\leqslant n.$ And we denote the determinant of the Vandermonde matrix $[a_i^{j-1}]_{1\leqslant i,j\leqslant r}$ by $V(a_1,a_2,\cdots,a_r).$ In view of \eqref{pi}, we obtain
\begin{equation}
    \begin{aligned}
    \det N&= \prod_{k=1}^{n}\eta_k\cdot V(\eta_1,\eta_2,\cdots,\eta_n)\\
    &=\prod_{1\leqslant j<k\leqslant n }(\zeta^{k^2}-\zeta^{j^2}).
\end{aligned}
\end{equation}
Clearly, $M$ is the $(n+1,2)$-cofactor of the matrix
\begin{equation}\label{dn2}
   M(x)=
  \begin{bmatrix}
    1 & 1  & \cdots & 1 &1\\
    \eta_1 & \eta_2 & \cdots & \eta_n & x \\
 \eta_1^2 & \eta_2^2 & \cdots & \eta_n^2 & x^2 \\
 \vdots & \vdots & \ddots & \vdots &\vdots \\
 \eta_1^n & \eta_2^n & \cdots & \eta_n^n & x^n \\
  \end{bmatrix}.
\end{equation}
Observe that
 \begin{equation}
     \begin{aligned}
    \det M(x)&= V(\eta_1,\eta_2,\cdots,\eta_n,x) \\
    &=  \prod_{i=1}^n(x-\eta_i)\prod_{1\leqslant j<k\leqslant n}(\eta_k-\eta_j).
 \end{aligned}
 \end{equation}
 Expanding $\det M(x)$ along the last column, we find that the coefficient of $x$ in $\det M(x) $ is
$(-1)^{n+3} \det M.
$ By (3.5), we have
\begin{equation}
  \begin{aligned}
   \det M&= (-1)^{n+3}\prod_{1\leqslant j<k\leqslant n}(\eta_k-\eta_j)
  \cdot (-1)^{n-1}\prod_{i=1}^{n}\eta_i\sum_{l=1}^{n}\eta_l^{-1} \\
   &=\prod_{1\leqslant j<k\leqslant n}(\eta_k-\eta_j)\sum_{l=1}^{n}\eta_l^{-1}.
\end{aligned}
\end{equation}
Since
\begin{equation}\label{eta-1}
\sum_{l=1}^{n}\eta_l^{-1}=\overline{\sum_{l=1}^{n}\eta_l}=   (\overline{\tau}-1)/2,
\end{equation}
we have
\begin{equation}\label{n1}
    \det M=\frac{1}{2}(\overline{\tau}-1)\times \prod_{1\leqslant j<k\leqslant n }(\zeta^{k^2}-\zeta^{j^2}).
\end{equation}
Thus,
\begin{equation}\label{n1n2}
    \det M \cdot \det N=\frac{1}{2}(\overline{\tau}-1)\times \prod_{1\leqslant j<k\leqslant n }(\zeta^{j^2}-\zeta^{k^2})^2.
\end{equation}

Next, we evaluate $\textbf{v}^T\mathrm{adj}(N)^T\mathrm{adj}(M)\textbf{v}.$ Write $\mathrm{adj}(M)=[M_{ij}]_{1\leqslant i,j \leqslant n}$ and $\mathrm{adj}(N)^T=[N_{ij}]_{1\leqslant i,j \leqslant n}.$ Then
$$\mathrm{adj}(N)^T\mathrm{adj}(M)=[n_{ij}]_{1\leqslant i,j\leqslant n},
 $$
where
$$ n_{ij}=\sum_{k=1}^n N_{ik}\cdot M_{kj}.$$
Hence,
\begin{equation}
 \begin{aligned}
    \textbf{v}^T\mathrm{adj}(N)^T\mathrm{adj}(M)\textbf{v}&=\sum_{i=1}^n\sum_{j=1}^n n_{ij}  \\
    &=\sum_{i=1}^n\sum_{j=1}^n\sum_{k=1}^{n} N_{ik}\cdot M_{kj}  \\
    &=\sum_{k=1}^n \(\sum_{i=1}^nN_{ik}\)\(\sum_{j=1}^nM_{kj}\).
\end{aligned}
\end{equation}

For any matrix $A,$  let $A^{(k)}$ denote the matrix obtaining from $A$ via replacing all entries in the $k$-th column by $1$. Expanding $\det M^{(k)}$ and $\det N^{(k)}$ along the $k$-th column, we obtain that
$$\det M^{(k)}=\sum_{j=1}^nM_{kj}
\ \ \t{and}\ \ \det N^{(k)}=\sum_{j=1}^nN_{ik}.
$$
Therefore,
\begin{equation}\label{n1n2'}
    \textbf{v}^T\mathrm{adj}(N)^T\mathrm{adj}(M)\textbf{v}=\sum_{k=1}^n\det M^{(k)}\det N^{(k)}.
\end{equation}
Since
$$ N^{(k)}=\begin{bmatrix}
    \eta_1 & \eta_2 & \cdots & \eta_{k-1} & 1 &\eta_{k+1} & \cdots & \eta_n  \\
 \eta_1^2 & \eta_2^2 & \cdots & \eta_{k-1}^2 & 1 &\eta_{k+1}^2 & \cdots & \eta_n^2  \\
 \vdots &\vdots &\ddots &\vdots & \vdots &\vdots & \ddots & \vdots \\
 \eta_1^{n-1} & \eta_2^{n-1} & \cdots & \eta_{k-1}^{n-1} & 1 &\eta_{k+1}^{n-1} & \cdots & \eta_n^{n-1}  \\
 \eta_1^n & \eta_2^n & \cdots & \eta_{k-1}^n & 1 &\eta_{k+1}^n & \cdots & \eta_n^n  \\
\end{bmatrix},
$$
we obtain
\begin{equation}
    \begin{aligned}
 \det N^{(k)}&= \frac{1}{\eta_k} \prod_{i=1}^n\eta_i \cdot V(\eta_1,\eta_2,\cdots,\eta_{k-1},1,\eta_{k+1},\cdots,\eta_n)\\
 &=\zeta^{-k^2}\cdot \prod_{1\leqslant i<j\leqslant n}({\zeta^{j^2}-\zeta^{i^2}})\prod_{\substack{1\leqslant j \leqslant n\\j\neq k}} \frac{1-\zeta^{j^2}}{\zeta^{k^2}-\zeta^{j^2}}.
\end{aligned}
\end{equation}
Clearly,
$$ M^{(k)}=\begin{bmatrix}
    1 & 1 & \cdots & 1 & 1 &1 & \cdots & 1  \\
 \eta_1^2 & \eta_2^2 & \cdots & \eta_{k-1}^2 & 1 &\eta_{k+1}^2 & \cdots & \eta_n^2  \\
 \eta_1^3 & \eta_2^3 & \cdots & \eta_{k-1}^3 & 1 &\eta_{k+1}^3 & \cdots & \eta_n^3  \\
 \vdots &\vdots &\ddots &\vdots & \vdots &\vdots & \ddots & \vdots \\
 \eta_1^{n-1} & \eta_2^{n-1} & \cdots & \eta_{k-1}^{n-1} & 1 &\eta_{k+1}^{n-1} & \cdots & \eta_n^{n-1}  \\
 \eta_1^n & \eta_2^n & \cdots & \eta_{k-1}^n & 1 &\eta_{k+1}^n & \cdots & \eta_n^n  \\
\end{bmatrix}
$$
is the $(n+1,2)$-cofactor of the matrix
\begin{equation}\label{n1'x}
  M^{(k)}(x)=\begin{bmatrix}
    1 & 1 & \cdots & 1 & 1 &1 & \cdots & 1&1  \\
     \eta_1 & \eta_2 & \cdots & \eta_{k-1} & 1 &\eta_{k+1} & \cdots & \eta_n&x  \\
 \eta_1^2 & \eta_2^2 & \cdots & \eta_{k-1}^2 & 1 &\eta_{k+1}^2 & \cdots & \eta_n^2&x^2  \\
 \vdots &\vdots &\ddots &\vdots & \vdots &\vdots & \ddots & \vdots&\vdots \\
 \eta_1^{n-1} & \eta_2^{n-1} & \cdots & \eta_{k-1}^{n-1} & 1 &\eta_{k+1}^{n-1} & \cdots & \eta_n^{n-1} &x^{n-1} \\
 \eta_1^n & \eta_2^n & \cdots & \eta_{k-1}^n & 1 &\eta_{k+1}^n & \cdots & \eta_n^n &x^n \\
\end{bmatrix}.
\end{equation}
Observe that
\begin{equation}
    \begin{aligned}
    \det M^{(k)}(x)&=V(\eta_1,\eta_2,\cdots,\eta_{k-1},1,\eta_{k+1},\cdots,\eta_n,x)\\
    &=(x-1)\prod_{\substack{j=1\\j\neq k}}^n(x-\zeta^{j^2}) \prod_{1\leqslant i<j\leqslant n}({\zeta^{j^2}-\zeta^{i^2}})\prod_{\substack{1\leqslant j \leqslant n\\j\neq k}} \frac{1-\zeta^{j^2}}{\zeta^{k^2}-\zeta^{j^2}} ,
\end{aligned}
\end{equation}
 and the coefficient of $x$ in $\det M^{(k)}(x)$ is
$ (-1)^{n+3}\det M^{(k)}.$
By \eqref{eta-1} and (3.14), we obtain that
    \begin{align*}
 \det M^{(k)}=&  (-1)^{n+3}(-1)^{n-1}\prod_{i=1}^n\eta_i\cdot \(\sum_{\substack{j=1\\j\neq k}}^n\frac{1}{\eta_k\eta_j}+\frac{1}{\eta_k}\) \\
 &\times\prod_{1\leqslant i<j\leqslant n}({\eta_j-\eta_i})\prod_{\substack{1\leqslant j \leqslant n\\j\neq k}} \frac{1-\eta_j}{\eta_k-\eta_j}  \\
 =& \frac{1}{\eta_k}\(1-\frac{1}{\eta_{k}}+\sum_{i=1}^n\eta_i^{-1}\) \times\prod_{1\leqslant i<j\leqslant n}({\eta_j-\eta_i})\prod_{\substack{1\leqslant j \leqslant n\\j\neq k}} \frac{1-\eta_j}{\eta_k-\eta_j} \\
 =&\frac{1}{\zeta^{k^2}}\((\overline{\tau}+1)/2-\frac{1}{\zeta^{k^2}}\)\prod_{1\leqslant i<j\leqslant n}({\zeta^{j^2}-\zeta^{i^2}})\prod_{\substack{1\leqslant j \leqslant n\\j\neq k}} \frac{1-\zeta^{j^2}}{\zeta^{k^2}-\zeta^{j^2}}.
\end{align*}
Thus, from (3.12) and the above equality,  we have
\begin{equation}
    \begin{aligned}
   &\textbf{v}^T\mathrm{adj}(N)^T\mathrm{adj}(M)\textbf{v}  \\
   =&\prod_{1\leqslant i<j\leqslant n}({\zeta^{i^2}-\zeta^{j^2}})^2\sum_{k=1}^n\frac{1}{\zeta^{2k^2}}\((\overline{\tau}+1)/2-\frac{1}{\zeta^{k^2}}\)\prod_{\substack{1\leqslant j \leqslant n\\j\neq k}} \(\frac{1-\zeta^{j^2}}{\zeta^{k^2}-\zeta^{j^2}}\)^2.
\end{aligned}
\end{equation}
Combining   \eqref{bp'}, (3.2), \eqref{n1n2} and (3.16), we obtain that
\begin{equation}
  \begin{aligned}
  \det A_p(x)=&  \(\overline{\tau}-1+
    (x\tau+1)\sum_{k=1}^n\(\frac{(\overline{\tau}+1)/2}{\zeta^{2k^2}}-\frac{1}{\zeta^{3k^2}}\)\prod_{\substack{1\leqslant j \leqslant n\\j\neq k}} \(\frac{1-\zeta^{j^2}}{\zeta^{k^2}-\zeta^{j^2}}\)^2\) \\
    &\times \frac{2^{n-1}}{\tau^{n}}\prod_{1\leqslant j<k\leqslant n}({\zeta^{j^2}-\zeta^{k^2}})^2.
\end{aligned}
\end{equation}

Now we distinguish two cases.

{\it Case} 1. $p\equiv 1 \pmod 4$.

In this case, by  Lemma 2.4 we have
 \begin{equation}
     \begin{aligned}
   &\sum_{k=1}^n\(\frac{(\overline{\tau}+1)/2}{\zeta^{2k^2}}-\frac{1}{\zeta^{3k^2}}\)\prod_{\substack{1\leqslant j \leqslant n\\j\neq k}} \(\frac{1-\zeta^{j^2}}{\zeta^{k^2}-\zeta^{j^2}}\)^2   \\
   =&\frac{\sqrt{p}+1}{2}\(\varepsilon_p^{-2h(p)}-1\)-(\sqrt{p}+1)\varepsilon_p^{-2h(p)}+1  \\
   =&-\frac{\sqrt{p}+1}{2}\varepsilon_p^{-2h(p)}-\frac{\sqrt{p}-1}{2}.
 \end{aligned}
 \end{equation}
 Combining Lemma 2.2, (3.17) amd (3.18) we have
 \begin{equation}
      \begin{aligned}
   \det  A_p(x)=&\(\sqrt{p}-1+
   (\sqrt{p}x+1)\(-\frac{\sqrt{p}+1}{2}\varepsilon_p^{-2h(p)}-\frac{\sqrt{p}-1}{2}\) \)\\
   &\times \frac{2^{n-1}}{\sqrt{p}^{n}} (-1)^{(p-1)/4}p^{(p-3)/4}\varepsilon_p^{h(p)}  \\
   =& \((\sqrt{p}-1)\varepsilon_p^{h(p)}-
   (\sqrt{p}x+1)\(\frac{\sqrt{p}+1}{2}\varepsilon_p^{-h(p)}+\frac{\sqrt{p}-1}{2}\varepsilon_p^{h(p)} \)  \)\\
   &\times \frac{2^{(p-3)/2}}{\sqrt{p}}\(\frac{2}{p}\)
   .
 \end{aligned}
 \end{equation}

 As $\varepsilon_p$ has norm 1 and $h(p)$ is odd, from
 \cite[Chapter 6, Theorems 4 and 6]{Co} we obtain
 \begin{equation}\label{-hp}
     \varepsilon_p^{-h(p)}=-a_p+\sqrt{p}b_p.
 \end{equation}
In view of \eqref{vp} and \eqref{-hp}, it is easy to verify that
\begin{equation}
    \begin{aligned}
   & \frac{\sqrt{p}+1}{2}\varepsilon_p^{-h(p)}+\frac{\sqrt{p}-1}{2}\varepsilon_p^{h(p)}\\
   =&\frac{\sqrt{p}+1}{2}(-a_p+\sqrt{p}b_p)+\frac{\sqrt{p}-1}{2}(a_p+\sqrt{p}b_p)  \\
    =& -a_p+pb_p.
\end{aligned}
\end{equation}
Combining \eqref{vp}, (3.19) and the above equality, we obtain
\begin{align*}
    \det  A_p(x)=&\frac{2^{(p-3)/2}}{\sqrt{p}}\(\frac{2}{p}\)\times\((\sqrt{p}-1)(a_p+\sqrt{p}b_p)-(\sqrt{p}x+1)(-a_p+pb_p)\) \nonumber\\
    =&\frac{2^{(p-3)/2}}{\sqrt{p}}\(\frac{2}{p}\)\times\(\sqrt{p}(a_p-pb_p)x+\sqrt{p}(a_p-b_p)\) \nonumber\\
    =&  2^{(p-3)/2}\(\frac{2}{p}\)(a_p-b_p+(a_p-pb_p)x).
\end{align*}

{\it Case} 2. $p\equiv 3 \pmod 4$.

In this case, by Lemma 2.5, we obtain
\begin{equation}
    \begin{aligned}
   &\sum_{k=1}^n\(\frac{(\overline{\tau}+1)/2}{\zeta^{2k^2}}-\frac{1}{\zeta^{3k^2}}\)\prod_{\substack{1\leqslant j \leqslant n\\j\neq k}} \(\frac{1-\zeta^{j^2}}{\zeta^{k^2}-\zeta^{j^2}}\)^2   \\
   =&\frac{1}{2}(-i\sqrt{p}+1)(-2)-(-2+i\sqrt{p}) =1
\end{aligned}
\end{equation}
Combining Lemma 2.2, (3.17) and (3.22), we have
\begin{align}
   \det  A_p(x)=& \frac{2^{(p-3)/2}}{(i\sqrt{p})^{(p-1)/2}}(-p)^{(p-3)/4}(-i\sqrt{p}-1+i\sqrt{p}x+1) \nonumber \\
   =& 2^{(p-3)/2}(x-1). \nonumber
\end{align}
This ends the proof.
\qed

\setcounter{conjecture}{0}
\end{document}